\documentclass[a4paper, 10pt, oneside]{article}

\title{Number variance from a probabilistic perspective:
infinite systems of independent Brownian motions and symmetric $\alpha$-stable
processes.}

\author{Ben Hambly\footnote{E-Mail: hambly@maths.ox.ac.uk} \,and Liza Jones\footnote{E-Mail: jones@maths.ox.ac.uk. Research supported by the EPSRC through the Doctoral Training Account scheme} \\  \small{Mathematical Institute, University of Oxford.}}

\date{July 2006}

\usepackage{amsmath}
\usepackage{amsfonts}
\usepackage{mathrsfs}
\usepackage{amsthm}
\usepackage{epsfig}

\numberwithin{equation}{section}

\newcommand{\ud}{\mathrm{d}}

\newtheorem{lem}{Lemma}
\newtheorem{prop}{Propostition}
\newtheorem{cor}{Corollary}

\newtheorem{thm}{Theorem}

\newtheorem{rem}{Remark}

\begin{document}

\maketitle

\begin{abstract}

Some probabilistic aspects of the number variance statistic are
  investigated.  Infinite systems of independent Brownian motions and
  symmetric $\alpha$-stable processes are used to construct new
  examples of processes which exhibit both divergent and saturating
  number variance behaviour. We derive a general expression for the
  number variance for the spatial particle configurations arising from
  these systems and this enables us to deduce various limiting
  distribution results for the fluctuations of the associated counting
  functions. In particular, knowledge of the number variance allows us
  to introduce and characterize a novel family of centered, long
  memory Gaussian processes. We obtain fractional Brownian motion as a weak limit of
  these constructed processes.
\end{abstract}


\section*{Introduction}

Let $(X, \mathcal{F}, \mathbb{P})$ be a point process on $\mathbb{R}$,
that is, a collection $X:=\{(x_i)_{i=-\infty}^{\infty} \,:\, x_i \in
\mathbb{R},\,\, \forall i \,\,\textrm{ and }\, \#(x_i : x_i \in [b,
b+L]) < \infty \,\,\forall \, b\in \mathbb{R}, L\in\mathbb{R}_+\}$, with
$\mathcal{F}$ the minimal $\sigma$-algebra generated by these point
configurations and $\mathbb{P}$ some probability measure on $(X,
\mathcal{F})$.  The associated $\textit{number variance}$ is defined
as
\begin{equation*}
Var(N[b, b+L]):=\textrm{ Variance}_{\mathbb{P}}[\# (x_i : x_i \in [b,b+L])].
\end{equation*}
More generally, in order to deal with non-spatially homogeneous cases,
it is more convenient to work with the averaged number variance which
we define as
\begin{eqnarray}
V[L] := \mathbb{E}_b[Var(N[b,b+L])] \nonumber
\end{eqnarray}
taking an appropriate uniform average of the number variance over different
intervals of the same length.
As we increase the length
$L$ of the underlying interval the number of points
contained in that interval will also increase. However, in many
situations, it is not immediately clear what will happen to the
variance of this number as the interval length grows. One of
the main questions considered in this paper will be the behaviour of
$V[L]$ as $L \rightarrow \infty$.
We shall see that, somewhat counter-intuitively, in some instances we have
\begin{equation*}
\lim_{L \rightarrow \infty} V[L]= \kappa < \infty,
\end{equation*}
in which case we will say that the number variance \textit{saturates}
to the level $\kappa \in \mathbb{R}_+$.

Until recently, interest in the number variance statistic had been
confined to the fields of random matrix theory (see e.g.
\cite{mehta,guhr2}) and quantum theory (e.g \cite{berry3, sarnak})
where it is commonly used as an indicator of spectral rigidity.  In
the study of quantum spectra, the manner in which the number variance
(of eigenvalues) grows with interval length provides an indication of
whether the underlying classical trajectories are integrable, chaotic
or mixed. In the large energy limit, the spectral statistics of a
quantum system with strongly chaotic classical dynamics, should,
according to the conjecture of Bohigas, Giannoni and Schmidt
\cite{bogschmit}, agree with the corresponding random matrix
statistics which are independent of the specific system. However, in
reality, for many such systems, when it comes to the long range
``global" statistics, this random matrix universality breaks down. In
these cases, the number variance typifies this transition to
non-conformity in that following an initial random matrix consistent
logarithmic growth with interval length, the number variance then
saturates.

Attempts to improve the understanding of the deviations from random
matrix predictions have led to convincing explanations for number
variance saturation behaviour in terms of periodic orbit theory, see
for example \cite{berry, berrykeating, aurichsteiner}. In \cite{berry}
there is a heuristic derivation of an explicit formula for the
empirical number
variance of the zeros of the Riemann zeta function on the critical
line which is consistent with numerical evidence.


In the last few years, number variance has been considered from a
slightly different viewpoint in relation to determinantal point
processes (see \cite{soshnikov} for background on this topic).
Results on the growth of the number variance for these processes are
given in e.g. \cite{soshnikov,sosh:airybess}. The emphasis in these
cases is on ascertaining the divergence of the number variance as
this, it turns out, is the key ingredient needed to prove Gaussian
fluctuation results for the counting functions of a large class of
determinantal processes, including those arising naturally in the
context of random matrix theory.

Motivated by the Riemann zeta example, Johansson \cite{kurt:varsat}
recently introduced an example of a point process with
determinantal structure which demonstrates the same type of number
variance saturation behaviour as conjectured by Berry (\cite{berry})
for the Riemann zeroes.  This process is constructed from a system of
$n$ non-colliding Brownian particles started from equidistant points
$u_j= \Upsilon + a(n-j)$ with $\Upsilon \in \mathbb{R},\quad a \in
\mathbb{R}_+ \quad j=1, \ldots, n.$ There are a number of approaches
to describing such a system, see \cite{kurt:wigner, joc} for details.
In any case, it can be shown that the configuration of particles in
space formed by the process at a fixed time $t$, is a determinantal
process and as such its correlation functions (or joint intensities)
take the form
\begin{align}
R^{(n)}_m(x_1,x_2, \ldots, x_m)= \textrm{det}\Big(K_t^{(n)}(x_i,x_j)\Big)_{i,j =1}^m.
\end{align}
Here and for determinantal processes in general, the determinant may
be interpreted as the density of the probability of finding $m$ of the
points in the infinitesimal intervals around $x_1,x_2, \ldots, x_m$.

As the number of particles $n \rightarrow \infty$, Johansson shows
that the correlation kernel $K_t^{(n)}$ converges uniformly to a
kernel with leading term (the correction is of order $e^{-d \beta},
\beta >0$)
\begin{align}
K_t( ax, ay) = \frac{1}{a} \,\Big{(}
\frac{\sin \pi (x-y)}{\pi (x-y)} + \frac{d \cos \pi (x+y) + (y-x)
\sin \pi (y+x)}{\pi( d^2 + (y-x)^2)} \Big{)}, \label{leadingpart}
\end{align}
where $d= \frac{2 \pi t}{a^2}$.  This limiting kernel defines, at each
time $t$, a limiting, non-spatially homogeneous determinantal process
which may loosely be thought of as the configuration of particles in
space formed by an infinite system of non-colliding Brownian particles
started from an equispaced initial configuration.  When the interval
length $L$ is small relative to $d$, the averaged number
variance for this process computed from the leading part of the
correlation kernel \eqref{leadingpart} has leading term
\begin{eqnarray}
\frac{1}{\pi^2} (\log (2 \pi L/a) +
\gamma_{\textrm{\tiny{Euler}}}+1). \label{varsine}
\end{eqnarray}
 However, if $d$
is held constant, while $L$ is increased, it is deduced that the
number variance saturates to the level
\begin{eqnarray}
\frac{1}{\pi^2} ( \log(2\pi d) + \gamma_{\textrm{\tiny{Euler}}} +1).
\end{eqnarray}
The ``small L" expression \eqref{varsine} agrees with the number
variance of the determinantal point process associated with the sine
kernel of density $a$
\begin{align}
K(x_i,x_j) = \frac{\sin \pi (x_i - x_j)/a}{ \pi(x_i - x_j)}. \label{sinekernel}
\end{align}
This determinantal process is the universal scaling limit obtained for
the eigenvalues of random matrices from e.g. the Gaussian Unitary
Ensemble and U(n), as matrix size tends to infinity (see for example,
\cite{mehta}).

For our purposes it will be convenient to think of the above averaged
number variance as the number variance of the spatial particle
configurations arising from an ``averaged model" which we choose to
interpret as an infinite system of non-colliding Brownian motions
started from the initial positions
\begin{eqnarray}
u_j= a(j-\epsilon), \qquad j \in \mathbb{Z}, \quad a \in \mathbb{R}_+,
\quad \epsilon \sim \textrm{Uniform}[0,1]. \label{initialconfig}
\end{eqnarray}

The two key features of Johansson's model as just described are the equispaced
starting positions and the fact that the Brownian particles are
conditioned to not collide.  In this work we consider
the number variance statistic in the independent
process analogues of Johansson's model. Since these independent
versions do not have determinantal structure, existing number variance
results no longer apply.

The paper is organized as follows. We begin by deriving a general
formula for the number variance for the spatial particle
configurations arising from
infinite systems of independent Brownian motions and symmetric
$\alpha$-stable processes in $\mathbb{R}$ started from the initial
configuration \eqref{initialconfig}. This enables us to
deduce the asymptotic behaviour of the statistic as the length of the interval over
which it is defined goes to infinity. We give an explicit formula for
the saturation level for the cases in which saturation does occur.
Once this is achieved we are able to explain the number variance
saturation phenomenon in terms of the tail distribution behaviours of
the underlying processes. We provide two specific illustrative
examples as corollaries and briefly discuss how our results compare to
those already existing in the literature. We conclude the first
section by demonstrating the close relationship between the number
variance and the rate of convergence of the distribution of the
associated counting function to a Poisson law.

In the second section we use the
number variance to prove Gaussian fluctuation results for the counting
functions of our particle configurations in two different scalings.
In the third and final section we add some dynamics to the
fluctuations of the counting functions to construct a collection of
processes, each of which is shown to converge weakly in $C[0,\infty)$
to a centered Gaussian process with covariance structure similar in
form to that of a fractional Brownian motion. Our earlier results on the
behaviour of the number variance allow us to better characterize these
limiting processes. In particular, the long-range dependence property
exhibited by the covariance of their increments is directly determined
by the rate of growth of the associated number variance. In the cases
corresponding to $\alpha \in (0,1)$, a further rescaling of the limiting Gaussian
processes allow us to recover fractional Brownian motions of Hurst parameters
$\frac{1-\alpha}{2}$ as weak limits.

\section{The independent particle cases}

\subsection{A Poisson process initial configuration}

We begin by illustrating the effect of the initial positions on the
number variance of the spatial particle configurations arising from
such infinite systems of processes as those considered in this paper.
The following theorem is the well known result (see for example,
\cite{doob}, Chapter VIII, section 5) that the Poisson
process is invariant for an infinite system of independent particles
with the same evolution.
%

\begin{thm} \label{poissonthm}
Consider an infinite collection of independent, identical in law, translation invariant
real-valued stochastic processes $\{(X_j(t), t \geq 0); j \in
\mathbb{Z}\}$.
Suppose that $\{X_j(0)\}_{j=-\infty}^{\infty}$ is a Poisson process of
intensity $\,\theta$ on $\mathbb{R}$. Then $\{X_j(t)\}_{j=-\infty}^{\infty}$ is a Poisson
process of intensity $\theta$ for every $t$.
\end{thm}

Consequently, if we begin with a ``mixed up" Poisson process initial
configuration and allow each particle to move as (say) a L\'evy
process, independently of the others then we observe a non-saturating,
linearly growing number variance ($V[L]= \theta L$) at the start and
for all time.

\subsection{The Brownian and symmetric $\alpha$-stable cases}

This last theorem served to highlight the importance of the regularity
or rigidity of the starting configuration \eqref{initialconfig} in
determining the number variance behaviour.  However, it is reasonable
to suppose, that in Johansson's model, the strong restrictions placed
on the movement of the Brownian particles must also contribute
significantly to the saturation behaviour that its number variance
demonstrates.
This leads us to ask; what would happen
if we started with an initial configuration such as
\eqref{initialconfig} but did not place such strong restrictions on
the movement of the particles?  We answer this question for the cases
in which each particle moves independently as a one-dimensional
Brownian motion or as a symmetric $\alpha$-stable process on
$\mathbb{R}$.

Recall (see e.g. \cite{sato})
that a symmetric $\alpha$-stable process is a L\'evy process
$(X^{\alpha,c}(t),t\geq 0)$ with characteristic function, for each
$t \geq 0$, given by
\begin{eqnarray}
\phi_{X^{\alpha,c}(t)}(\theta)&:=&\mathbb{E}\big[e^{i \theta
X^{\alpha,c}(t)}\big]\nonumber
\\&=&\exp\big(-t\, c \, |\theta|^{\alpha}\, \big), \qquad c>0,
\quad \alpha \in (0,2). \label{charfunc}
\end{eqnarray}
Some of the properties enjoyed by this class of processes are;
\begin{itemize}
\item $\{X^{\alpha,c}(t), t \geq 0\}$ with associated transition
  density $p_t(x,y) \,,\, x,y \in \mathbb{R}$ is temporally and spatially homogeneous.
\item $\{X^{\alpha,c}(t)\}$ is symmetric and self-similar in the sense that
$\{X^{\alpha,c}(t)\}\,\, \buildrel {\textrm{dist}} \over {=} \,\, \{-X^{\alpha,c}(t)\}$ and
$\{\lambda^{-1/\alpha} X^{\alpha,c}\big(\lambda t\big)\}\,\,\buildrel {\textrm{dist}} \over
{=}\,\, \{X^{\alpha,c}(t)\}$ for constants $\lambda \in \mathbb{R}$.
\end{itemize}
The arguments that follow also apply to the $\alpha =2$ Gaussian
cases. Note that we have standard Brownian motion when $\alpha=2, c=
\frac{1}{2}$.

\begin{thm} \label{sas}
Fix a symmetric $\alpha$-stable process $(X^{\alpha,c}(t), t \geq
0)$ on $\mathbb{R}$ with properties as described above. Suppose we
start an independent copy of this process from each of the starting
positions
\begin{eqnarray}
 u_j := a(j - \epsilon), \quad j \in \mathbb{Z}, \quad
 \textrm{where}\,\, a \in \mathbb{R}_+ \,\, \textrm{and} \,\,\, \epsilon \sim
 \textrm{Uniform}[0,1]. \label{startpos}
 \end{eqnarray}
 The configuration of particles in space formed by this infinite
 system of independent symmetric $\alpha$-stable processes at a fixed time $t$ has
 number variance
 \begin{align}
V_t^{\alpha,c,a}[L] = \frac{L}{a} + \frac{2}{\pi}\int_0^{\infty} \frac{e^{-2 c t
 (\theta/a)^{\alpha}}}{\theta^2}\Big[\cos\Big(\frac{L
 \theta}{a}\Big)-1\Big] \, \ud \theta.  \label{numvarsas}
 \end{align}
\end{thm}
\begin{proof}
  Let $\{(X_j^{\alpha,c}(t),t \geq 0), j \in \mathbb{Z} \}$ be the
  independent copies of the chosen symmetric $\alpha$-stable process
  indexed by $j \in \mathbb{Z}$.  Denote the law of each
  $X_j^{\alpha,c}$ started from $x \in \mathbb{R}$ by $\mathbb{P}^x$,
and write $\mathbb{P}:= \mathbb{P}^0$.
Now the number of symmetric $\alpha$-stable particles in an interval
$[0,L] \subset \mathbb{R}$ at time $t$ is given by the sum of
indicator random variables
\begin{equation}
N_t^{\alpha,c,a}[0,L] = \sum_{j=-\infty}^{\infty}
\mathbb{I}[X_j^{\alpha,c}(t)+u_j \in [0,L]], \label{sumofi}
\end{equation}
where $(u_j)_{j=-\infty}^{\infty}$ is given by (\ref{startpos}).
Note that by construction, for this ``averaged model", for all $b \in \mathbb{R}$
   we have
 $N_t^{\alpha, c,a}[0,L] \buildrel {\textrm{dist}} \over {=} N_t^{\alpha, c,a}[b,
 b+L]$. Thus the number variance is given by
\begin{eqnarray}
V_t^{\alpha,c,a}[L] &:=& \textrm{Var}\big[N_t^{\alpha,c,a}[0,L]
\big]\nonumber \\&=& \!\!\!\!\! \sum_{j=-\infty}^{\infty}
\!\!\mathbb{P}[X_j^{\alpha,c}(t)+u_j \in
[0,L]]\,\mathbb{P}[X_j^{\alpha,c}(t)+u_j \notin [0,L]].
\label{firstvar}
\end{eqnarray}
We can use the self-similarity property and the independence of
$\epsilon$ and the $X_j$ to write the probabilities under
consideration as convolutions which then allows us to deduce
\begin{eqnarray}
 \sum_{j=-\infty}^{\infty} \negthickspace \mathbb{P}[X_j^{\alpha,c}(t) + a(j -
\epsilon) \in [0,L]] \nonumber
=  \int_{-\infty}^{\infty} \int_0^{L/a} p_{t/a^{\alpha}}(x,y) \, \ud
y \, \ud x.
\nonumber
\end{eqnarray}
Hence
\begin{eqnarray}
&&V_t^{\alpha,c,a}[L] \label{eq:vtl} \\ &=&
\underbrace{\int_{-\infty}^{\infty}\! \int_0^{L/a} \!\!\!p_{t/a^{\alpha}}
(x,y) \, \ud y \, \ud x}_{T_1} -
\underbrace{\int_{-\infty}^{\infty}\!\int_0^{L/a}\!\! \!\int_0^{L/a}
\!\!\!p_{t/a^{\alpha}}(x,y)p_{t/a^{\alpha}}(x,z)\, \ud z \, \ud y \, \ud
x}_{T_2}. \nonumber
\end{eqnarray}
By Fubini's Theorem and symmetry we have
\begin{equation}
T_1 = \int_0^{L/a} \int_{-\infty}^{\infty} p_{t/a^{\alpha}}(y,x)
\, \ud x \, \ud y
= \frac{L}{a}.\label{eq:upbdvtl}
\end{equation}
For the other term we make use of the Chapman-Kolmogorov identity and
the spatial homogeneity before performing an integral switch to obtain
\begin{align}
T_2 &= \int_0^{L/a} \int_{-y}^{(L/a)-y}p_{2t/a^{\alpha}}(0,z) \, \ud z \,
\ud y
\nonumber
\\ &= \frac{L}{a} \, \int_{-L/a}^{L/a} p_{2t/a^{\alpha}}(0,z) \, \ud z
+ \int_{-L/a}^0 \!\!z \, p_{2t/a^{\alpha}}(0,z) \, \ud z - \int_0^{L/a}\!\!\!
z \, p_{2t/a^{\alpha}}(0,z) \, \ud z. \nonumber
\end{align}
From the symmetry property we know that, for each $t$,
$p_t(0,z)$ is an even function in $z$ and
$g(z) := z \, p_t(0,z)$ is an odd function in $z$. These facts allow
us to conclude, bringing the two terms together,

\begin{equation} \label{genform}
V_t^{\alpha,c,a}[L] = \frac{2 L}{a} \int_{L/a}^{\infty}p_{2t/a^{\alpha}}(0,z) \, \ud z
+ 2 \int^{L/a}_0 z\,p_{2t/a^{\alpha}}(0,z) \, \ud z.
\end{equation}
Applying Fourier inversion to the characteristic function
$\phi_{X^{\alpha,c}(t)}(\theta)$ given at (\ref{charfunc}), we deduce that the transition
density can be expressed as
\begin{equation*}
p_t(0,z | \alpha, c) = \frac{1}{\pi} \int_0^\infty \cos(z \theta)
e^{-ct\theta^{\alpha}} \, \ud \theta.
\end{equation*}
Using this density, the symmetry property and the expression \eqref{genform}
we obtain
\begin{eqnarray}
V_t^{\alpha,c,a}[L] = \frac{L}{a} + \frac{2}{\pi} \int_0^{\infty}
\frac{e^{-2ct(\theta/a)^{\alpha}}}{\theta^2}\Big[\cos\big(\frac{L\theta}{a}\big)
-1 \Big] \, \ud \theta, \nonumber
\end{eqnarray}
as required.
\end{proof}

Having found a general expression for the number variance we are able
to consider its behaviour as the interval length $L$ is increased.
Note that given positive real valued functions $g$ and $h$ we let $g
\sim h$ signify that $\lim \frac{g}{h} = 1$.

\begin{thm}\label{satlevel:thm}
Consider the number variance $V_t^{\alpha,c,a}[L]$ for the system of
symmetric $\alpha$-stable processes considered above. We claim that
\begin{equation} \label{asymvar}
V_t^{\alpha,c,a}[L] \sim \begin{cases} \quad k_{\alpha,c,a,t}\,
  L^{1-\alpha}, & \textrm{for }\alpha \in(0,1), \\
\quad k_{\alpha,c,a,t}\,\log(L), & \textrm{for } \alpha =1, \\
\quad k_{\alpha,c,a,t}\,L^{1-\alpha} \,+\, \kappa_{sat}(\alpha,c,a,t),
& \textrm{for }\alpha \in(1,2), \\
\quad k_{\alpha,c,a,t}\,e^{-L^2/8ct} \,+\, \kappa_{sat}(\alpha,c,a,t),
& \textrm{for } \alpha=2,
\end{cases}
\end{equation}
as $L \rightarrow \infty$, where $k_{\alpha,c,a,t}$ is a constant
depending on $\alpha,c,a$ and $t$, and
\begin{equation}
\kappa_{sat}(\alpha,c,a,t) = \frac{2}{a \pi}(2tc)^{1/\alpha}
\Gamma\Big(1-\frac{1}{\alpha}\Big),\label{satlevelsas}
\end{equation}
where $\Gamma(x) := \int_0^{\infty} s^{x-1} \, e^{-s} \, \ud s$ is
the usual Gamma function.
\end{thm}
\begin{proof} From (\ref{genform}) note that we may
re-write the expression for the number variance as
\begin{align}
V_t^{\alpha,c,a}[L] = \frac{L}{a}\, \mathbb{P}
\Big[\Big|X^{\alpha,c}\Big(\frac{2t}{a^{\alpha}}\!\Big)\Big|\! >\! L/a \Big] \! + \!
\mathbb{E}_{\mathbb{P}}\Big[\Big|X^{\alpha,c}\Big(\frac{2t}{a^{\alpha}}\Big)\Big| \,
; \, \Big|X^{\alpha,c}\Big(\frac{2t}{a^{\alpha}}\!\Big)\Big| \!< \!L/a \Big],
\label{numvardiffform}
 \end{align}
Now the behaviour of the first term in this expression is well known (see \cite{samtaqqu}, page 16). For $\alpha \in (0,2)$ we have
\begin{align*}
\frac{L}{a} \,\mathbb{P} \Big[\Big|X^{\alpha,c}\Big(\frac{2t}{a^{\alpha}}\Big)\Big|
> L/a \Big] \sim k_{\alpha}\frac{2ct}{a}\,
L^{1-\alpha} \quad \textrm{as} \,\,\,L
\rightarrow \infty,
\end{align*}
where $k_{\alpha}= \Big(\int_0^{\infty} x^{-\alpha} \sin x \, \ud
x\Big)^{-1}$. We use $\log(L)$ when $\alpha =1$. In the Gaussian case
($\alpha=2$) we have instead
\begin{eqnarray}
\frac{L}{a} \,\mathbb{P} \Big[\Big|X\Big(\frac{2t}{a^2}\Big)\Big| > L/a \Big] \sim
\frac{a}{\sqrt{2 c \pi t}}\, e^{-L^2/8ct} \quad \textrm{as} \quad L
\rightarrow \infty. \nonumber
\end{eqnarray}
To deal with the
second term of (\ref{numvardiffform}) observe that we have
\begin{align*}
\mathbb{E}_{\mathbb{P}}\Big[\Big|X^{\alpha,c}(2t/a^{\alpha})\Big| \,
; \, \Big|X^{\alpha,c}(2t/a^{\alpha})\Big| < L/a \Big] \,\,  \buildrel L \rightarrow \infty \over  \rightarrow \,\,
&\mathbb{E}_{\mathbb{P}}\Big[\Big|X^{\alpha,c}(2t/a^{\alpha})\Big|\Big]\quad \\
= \,\, \,\, \,\,&
\int_0^{\infty}
\mathbb{P}\Big[\Big|X^{\alpha,c}(2t/a^{\alpha})\Big|>\lambda\Big] \ud
\lambda.
\end{align*}
Thus it is clear that the rate of divergence or saturation of the number variance is
determined by the upper tail of the distribution of the underlying
symmetric $\alpha$-stable process. This may also be seen by
differentiating the integral expression for $V_t^{\alpha,c,a}[L]$ with
respect to L
 \begin{align*}
 \frac{\partial V_t^{\alpha,c,a}[L]}{\partial L}
 = \frac{1}{a} \Big[\mathbb{P}\Big|X
 \big(\frac{2t}{a^{\alpha}}\big)\Big| \geq \frac{L}{a}\Big].
 \end{align*}
Consequently, for $\alpha \in (1,2]$ (the saturating cases)
\begin{align*}
\lim_{L \rightarrow \infty} V_t^{\alpha,c,a}[L] =
\mathbb{E}_{\mathbb{P}}\Bigg[\Bigg|X^{\alpha,c}\Big(\frac{2t}{a^{\alpha}}\Big)
\Bigg| \Bigg] =: \kappa_{sat}(\alpha,c,a,t) < \infty.
\end{align*}
The exact expression for $\kappa_{sat}$ can be obtained from the
moments. By \cite{shanbhag},
if $X$ is a symmetric $\alpha$-stable random variable with $0 <\alpha \leq 2$ and scale
$\sigma$, then for $-1 < \delta < \alpha$ we have
\begin{align*}
\mathbb{E}\Big[ \big| X \big|^{\delta}
\Big]=\frac{\sigma^{\delta/\alpha}\, 2^{\delta}\,
\Gamma\big(\frac{1+\delta}{2} \big) \,
\Gamma\big(1-\delta/\alpha\big)}{\Gamma(1/2)\,\Gamma(1-\delta/2)}.
\end{align*}
Applying this theorem with $\delta=1,
\sigma = \frac{2ct}{a^\alpha}$ gives
(\ref{satlevelsas}).
To see how this fits in with the integral expression for
$V_t^{\alpha,c,a}[L]$, it may be verified that for $\alpha \in (1,2]$
\begin{align*}
\kappa_{sat}(\alpha,c,a,t) = -
\frac{2}{\pi}\int_0^{\infty}\frac{e^{-2ct(\theta/a)^{\alpha}}}{\theta^2}
\, \ud \theta.
\end{align*}
\end{proof}

The saturation result is now a consequence of Theorem~\ref{satlevel:thm}.

\begin{cor}
\begin{enumerate}
\item If $\alpha \in (0,1]$, at each time $t>0$,
$V_t^{\alpha,c,a}[L]$  \textbf{diverges} as $L \rightarrow \infty$.
\item If $\alpha \in (1,2]$, at each time $t<\infty$,
$V_t^{\alpha,c,a}[L]$  \textbf{saturates} to the level
$\kappa_{sat}(\alpha,c,a,t)$ as $L \rightarrow \infty$.
\end{enumerate}
\end{cor}

\begin{rem}{\rm
Even from the simplest equation \eqref{firstvar} it is clear that the
largest contributions to the number variance come from the activity at
the edges of the interval under consideration. Thus intuitively, the
fatter the tails of the distributions concerned the greater the number
of particles that may be in the vicinity of these edges making these
substantial contributions, the slower the decline in this number as
the interval length increases and consequently the slower the decay in
the growth of the number variance. }
\end{rem}

We now apply Theorems \ref{sas} and \ref{satlevel:thm} to two well known examples.
\begin{cor}[Brownian case]
Consider an infinite system of Brownian particles started from the initial configuration
$(u_j)_{j=-\infty}^{\infty}$ as given at \eqref{startpos}. The number
variance for this process is
\begin{equation*}
V_t^{2,\frac{1}{2},a}[L]= \frac{2}{a}\Bigg[L \, \Phi\Big(\frac{-L}{\sqrt{2t}}\Big) + \sqrt{
\frac{t}{\pi}} \Big(\, 1- e^{-L^2/4t} \,\Big) \Bigg],
\end{equation*}
where
\begin{equation}
\Phi(x):= \frac{1}{\sqrt{2 \pi}} \, \int_{-\infty}^x e^{-y^2/2} \,
\ud y. \nonumber
\end{equation}
As $L \rightarrow \infty$ this number variance saturates exponentially quickly to the level
\begin{equation*}
\frac{2}{a} \, \sqrt{\frac{t}{\pi}}.
\end{equation*}
\end{cor}


\begin{cor}[Symmetric Cauchy case]
Consider
an infinite system of symmetric Cauchy processes started from the initial configuration
$(u_j)_{j=-\infty}^{\infty}$ as given at \eqref{startpos}. The number
variance for this process is
\begin{eqnarray}
V_t^{1,1,a}[L]= \frac{L}{a}\Big[1-\frac{2}{\pi}
\arctan\big[\frac{L}{2t}\big] \Big] + \frac{2t}{a \pi}\Big[\log\Big(
1+ \Big(\frac{L}{2t}\Big)^2 \Big) \Big]. \label{cauchyvar}
\end{eqnarray}
For ``large"
(relative to 2t) L, we have
\begin{equation*}
V_t^{1,1,a}[L] \approx \frac{4t}{a \pi}
\log\Big(\frac{L}{2t}\Big)\quad \nonumber
\end{equation*}
and so the number variance diverges at a logarithmic rate as $L \rightarrow \infty$.
\end{cor}

\begin{rem}
Coincidentally, in the symmetric Cauchy case, if we set $a=1$, $t= \frac{1}{4 \pi}$ we have
\begin{eqnarray}
V_{\frac{1}{4 \pi}}^{1,1,1}[L] \approx \frac{1}{\pi^2}\log(2\pi L) \nonumber
\end{eqnarray}
and so we see similar number variance behaviour as in
the sine kernel case (\ref{varsine}).
\end{rem}


\begin{rem}{\rm
In general, the number variance for the sine kernel determinantal
process and for Johansson's determinantal process described in the
introduction, can be expressed in the form \eqref{genform} with
$p_{2t/a^{\alpha}}(0,z)$ replaced by the functions $\frac{sin^2 (\pi
  z)}{a \pi^2 z^2}$ and $\frac{sin^2 (\pi  z)}{a \pi^2 z^2} +
\frac{d^2 - z^2}{2a \pi^2(d^2 + z^2)^2}$ respectively. In the latter
case the saturation level is again given by the ``expectation" term. }
\end{rem}

\begin{rem}{\rm
The Cauchy ($\alpha=1$) case has the slowest growing non-saturating
number variance amongst all those considered here. Analogously, the
sine kernel determinantal process has the slowest growing number
variance amongst all translation invariant determinantal processes
whose kernels correspond to projections (i.e. the Fourier transform of
the kernel is an indicator) see \cite{soshnikov}. }
\end{rem}

\begin{rem}\label{extreme} {\rm
At the other extreme note that as $\alpha\to 0$ we recover Poisson
behaviour in that
\begin{equation*}
\lim_{\alpha \rightarrow 0} V_t^{\alpha,c,a}[L] = \frac{L}{a}(1-e^{-2ct}).
\end{equation*}}
\end{rem}

\subsection{A Poisson approximation for the counting function}
At the beginning of this section we recalled that a system of
independent processes (satisfying fairly general conditions) started
from a Poisson process on $\mathbb{R}$ remains in a Poisson process
configuration and hence demonstrates a number variance linear in L,
for all time. Now from (\ref{eq:vtl}) and (\ref{eq:upbdvtl})
we deduce that
\begin{align}
V_t^{\alpha,c,a}[L] \leq L/a,
\qquad \forall\,\, \alpha,c,a,t. \label{markovbound}
\end{align}
From the integral expression for the number variance
(\ref{numvarsas}), we see that for fixed $L$
\begin{align*}
V_t^{\alpha,c,a}[L] \rightarrow L/a \quad \textrm{as } t \rightarrow \infty,
\end{align*}
and for each $t$, as $L$ is decreased
\begin{align*}
V_t^{\alpha,c,a}[L] \buildrel {L \rightarrow 0} \over {\sim} L/a.
\end{align*}
So in both these limiting cases (as well as in the $\alpha \rightarrow
0$ case c.f. Remark \ref{extreme}) the maximal linear number variance
is attained. The balance of the parameters $L,a, \alpha$ and $t$,
encapsulated by the number variance, determines how ``far away" from
being Poisson the distribution of $N_t^{\alpha, c,a }[0,L]$ actually
is. Below we make this observation precise.

Recall that given a measurable space $(\Omega, \mathcal{F})$ we
define the total variation distance $d_{TV}(\cdot,\cdot)$ between two probability
measures $\mu_1, \mu_2$ on $\Omega$ by
\begin{equation}
d_{TV}(\mu_1,\mu_2) := \sup_{F \in \mathcal{F}} | \mu_1(F)-\mu_2(F)
|. \nonumber
\end{equation}

\begin{prop}\label{totalvar}
Let $\mathcal{L}(N^{\alpha,c,a}_t[0,L])$ denote the law of the random
variable $N_t^{\alpha,c,a}[0,L]$ defined at \eqref{sumofi}.
Let $\textrm{Po}(L/a)$ denote the law of a Poisson random variable
with mean $L/a$. Then, for $L \geq 1$,
\begin{align*}
\frac{L- V_t^{\alpha,c,a}[L]a}{32L}\! \leq \!
d_{TV} \!\big(\mathcal{L}(N_t^{\alpha,c,a}[0,L]), \textrm{Po}(L/a)\!
\big)\! \leq \!(1-e^{-\frac{L}{a}})\!\frac{L- V_t^{\alpha,c,a}[L]a}{L}.
\end{align*}
\end{prop}
\begin{proof}
  The result is an application of a theorem of Barbour and Hall
  \cite{barbhall}. Their theorem states that if $A:=\sum_j
  \mathbb{I}_j$ is the sum of independent indicator random variables
  indexed by $j$ and $q_j(L) = \mathbb{E}[\mathbb{I}_j]$, $\lambda =
  \sum_j q_j(L)$ then if we denote the law of A by $\mathcal{L}(A)$ we
  have
\begin{eqnarray}
\frac{1}{32} \min(1,1/\lambda) \,\sum_j q_j(L)^2 \leq d_{TV} \big(
\mathcal{L}(A), \textrm{Po}(L/a)\big) \leq
\frac{1-e^{-\lambda}}{\lambda} \,\sum_j q_j(L)^2.  \nonumber
\end{eqnarray}
In our specific case we have
$N_t^{\alpha,c,a}[0,L]$ as the sum of independent indicator random
variables given by (\ref{sumofi}), $\lambda= L/a$ and $\sum_j q_j(L)^2
= \frac{L}{a} - V_t^{\alpha,c,a}[L]$.
\end{proof}

\begin{rem}{\rm
For a fixed $t$, the Poisson approximation becomes less accurate as $L
\rightarrow \infty$. The greater the value of $\alpha$ the faster the
quality of the approximation deteriorates. For $\alpha >1$, due
to the fact that the number variance saturates, the approximation of
the law of $N_t^{\alpha,c,a}[0,L]$ by a Poisson distribution of mean
$L/a$ becomes very poor with the total variation distance clearly
bounded away from zero.}
\end{rem}

\section{Gaussian fluctuations of the counting function}
Thus far we have been concerned with the variance of the counting
function $N_t^{\alpha,c,a}[0,L]$ \eqref{sumofi}.
 Of course this variance is, by definition, a
description of the fluctuation of $N_t^{\alpha,c,a}[0,L]$ around its
mean $\frac{L}{a}$. In this section we will further characterize these
fluctuations.

\begin{prop}\label{fluc}
Let $N_t^{\alpha,c,a}[0,L]$,
$V_t^{\alpha,c,a}[L]$ denote, as usual,  the counting function and
number variance. For the cases with $\alpha \in (0,1]$ we have that
\begin{eqnarray}
\frac{N_t^{\alpha,c,a}[0,L]-L/a}{\sqrt{V_t^{\alpha,c,a}[L]}}
\end{eqnarray}
converges in distribution to a standard normal random
variable as $L \rightarrow \infty$.
\end{prop}

\begin{proof}
 Recall
that the cumulants $c_k, \, k \in \mathbb{N}$ of a real valued
random variable $Y$ are defined
by
\begin{eqnarray}
\log \mathbb{E}[\exp(i \theta Y)] = \sum_{k=1}^{\infty} c_k \frac{(i
\theta)^k}{k\, !}. \nonumber
\end{eqnarray}
Using the independence of the individual symmetric
$\alpha$-stable processes and then applying the Taylor expansion of
$\log(1+x)$ about zero, we have
\begin{align*}
\log \mathbb{E}_{\mathbb{P}}[\exp(i
\theta N_t^{\alpha,c,a}[0,L])] &= \sum_{j=-\infty}^{\infty} \log
\Big[\big( e^{i \theta} -1 \big) q_j(L) + 1 \Big] \\ &= \sum_{m=1}^{\infty}
\frac{(e^{i\theta}-1)^m}{m}
(-1)^{m+1}\Big(\sum_{j=-\infty}^{\infty}q_j(L)^m\Big),
\end{align*}
where $q_j(L):= \mathbb{P}[X_j^{\alpha,c}(t)+u_j \in [0,L]]$ and
$X_j^{\alpha,c}(t)$ denotes as usual the underlying symmetric $\alpha$-stable
process labelled by j with $u_j$ the corresponding starting position.
Hence, the cumulants of $N_t^{\alpha,c,a}[0,L]$ are given by
\begin{eqnarray}
c_k = \frac{\ud^k}{\ud \theta^k} \sum_{m=1}^{\infty}
\frac{(e^{i\theta}-1)^m}{m}
(-1)^{m+1}\Big(\sum_{j=-\infty}^{\infty}q_j(L)^m\Big) \Big|_{\theta=0}.
\label{cum}
\end{eqnarray}
It is straightforward to see that
\begin{align*}
c_1= L/a, \qquad \qquad  c_2=
\sum_{j=-\infty}^{\infty} q_j(L) - q_j(L)^2 = V_t^{\alpha,c,a}[L]
\end{align*}
 give the mean and number variance respectively. More generally from the
equation (\ref{cum}) it is possible to deduce the following
recursive relation
\begin{equation}
c_k =\sum_{m=2}^{k-1} \beta_{k,m} c_m +(-1)^k(k-1)!
\sum_{j=-\infty}^{\infty}q_j(L)-q_j(L)^k,  \label{recursive}
\end{equation}
where $\beta_{k,m}$ are constant, finite, combinatorial coefficients
which will not be needed here. Now let
\begin{equation*}
Y_t^{\alpha,c,a}:= \frac{N_t^{\alpha,c,a}[0,L] - L/a}{\sqrt{V_t^{\alpha,c,a}[L]}}.
\end{equation*}
It is easily deduced that the cumulants $\tilde{c_k},\,
k \in \mathbb{N}$ of $Y_t^{\alpha,c,a}$ are given by
\begin{eqnarray}
\tilde{c_1} &=& 0, \nonumber \\ \tilde c_k &=&
\frac{c_k}{(V_t^{\alpha,c,a}[L])^{k/2}}, \qquad \textrm{for $k \geq
2$}. \nonumber
\end{eqnarray}
To prove the Proposition it is sufficient to show that in the limit as $L
\rightarrow \infty$, the cumulants correspond those of a Gaussian
random variable. That is, we have $\tilde c_3=\tilde c_4=\tilde c_5
=\cdots =0$. Equivalently, we need to show
\begin{eqnarray}
c_k = o((V_t^{\alpha,c,a}[L])^{k/2})=o(c_2^{k/2}), \qquad \textrm{for
$k \geq 3$}. \nonumber
\end{eqnarray}
We use an induction argument.

Suppose that $c_m = o(c_2^{m/2})$ for $m=3,\ldots, k-1$. Assume,
without loss of generality, that
$k$ is even. We use the inequality
\begin{eqnarray}
q_j(L)-q_j(L)^k &=& \sum_{l=1}^{k-1}q_j(L)^l-q_j(L)^{l+1} \nonumber
\\ &\leq& (k-1)\,(q_j(L)-q_j(L)^2)\label{ineq1b}
\end{eqnarray}
in conjunction with the recursive relation for $c_k$ given by \eqref{recursive} to
deduce
\[  \sum_{m=2}^{k-1} \beta_{k,m} c_m \quad \quad \leq \quad  c_k  \quad \leq \quad
 (k-1)!(k-1)c_2 +\sum_{m=2}^{k-1} \beta_{k,m} c_m. \]
From our induction supposition this implies that
\begin{align}
o(c_2^{k/2}) \quad \leq \quad c_k \quad \leq \quad
o(c_2^{k/2}) + (k-1)!(k-1)\,c_2. \label{eqinduct}
\end{align}
However, from the results of the previous section, we know that, for
these cases with $\alpha \in(0,1]$, for $k \geq 3$,
$c_2^{\frac{k-2}{2}}= V_t^{\alpha,c,a}[L]^{\frac{k-2}{2}} \rightarrow
\infty$ as $L \rightarrow \infty$. Thus from \eqref{eqinduct} $c_k =
o(c_2^{k/2})$ also.
Now using the same arguments as for the inequality \eqref{ineq1b} we find
\begin{eqnarray}
 -\frac{1}{\sqrt{c_2}} \leq &\frac{c_3}{c_2^{3/2}}& \leq
\frac{1}{\sqrt{c_2}}. \nonumber
\end{eqnarray}
Thus
we have
$c_3=o(c_2^{3/2})$. By the induction argument we can deduce that
$\frac{c_k}{(c_2)^{k/2}} \rightarrow 0$ as $L \rightarrow \infty$
for all $k \geq 3$ which concludes the proof.
\end{proof}
\begin{rem}{\rm
The divergence of the number variance is relied upon in a similar way
to prove the analogous Gaussian fluctuation results for a large
class of determinantal processes, see \cite{costinlebowitz,
  sosh:airybess, soshnikov}. In all cases the overall structure of
the proof is the same. We note that the Proposition could also
have been proven by applying the Lindberg-Feller Central Limit
Theorem (see e.g. \cite{feller}). }
\end{rem}

Proposition \ref{fluc} applies to the cases with $\alpha \in(0,1]$.
The following convergence in distribution result applies to all cases
with $\alpha \in (0,2]$ and is obtained by allowing both interval
length and time $t$ tend to infinity together in an appropriate way.

\begin{prop}
For any fixed $s \in [0, \infty)$ we have that
\begin{eqnarray}
\frac{N_t^{\alpha,c,a}[0, s t^{1/\alpha}] - s t^{1/ \alpha} /a}{t^{1/2
    \alpha}}
\label{distlimit}
\end{eqnarray}
converges in distribution as $t \rightarrow \infty$, to a normal
random variable with zero mean and variance $f^{\alpha,c,a}(s)$, where
\begin{eqnarray}
f^{\alpha,c,a}(s) &:=&  V_1^{\alpha,c,a}[s] \nonumber  \\
&=& \frac{4 s}{a \pi}\int_0^{\infty}  \frac{\sin^2(u/2)}{u^2}\,\big( 1- e^{-2 c
    (u/s)^{\alpha}}\big) \, \ud u. \label{fdefn}
\end{eqnarray}
\end{prop}
\begin{proof}
Since $V_t^{\alpha,c,a}[s t^{1/\alpha}] \rightarrow \infty$ as $t
\rightarrow \infty$, a similar argument as for the proof of
Proposition \ref{fluc} allows us to conclude that
\begin{eqnarray}
\frac{N_t^{\alpha,c,a}[0, s t^{1/\alpha}] - s t^{1/\alpha}
  /a}{\sqrt{V_t^{\alpha,c,a}[ s t^{1/\alpha}]}}
\end{eqnarray}
converges in distribution as $t \rightarrow \infty$ to a standard
normal random variable.
Making the change of variable $u= s t^{1/\alpha} \theta/a$ in the integral
expression for the number variance \eqref{numvarsas}, using
$\int_0^{\infty} \sin^2(u)/u^2 du = \pi/2$ and a double angle formula
yields
\begin{align*}
\frac{V_t^{\alpha,c,a}[st^{1/\alpha}]}{t^{1/\alpha}} =\frac{4 s}{a
  \pi}\int_0^{\infty}  \frac{\sin^2 (u/2)}{u^2}\,\big( 1- e^{-2 c
    (u/s)^{\alpha}}\big) \, \ud u =: f^{\alpha,c,a}(s).
\end{align*}
Note that from \eqref{markovbound} we know $f^{\alpha,c,a}(s)< \infty$
for all $s < \infty$ and so the result follows from the scaling
property of the Gaussian distribution.
\end{proof}

\section{The fluctuation process}
We proceed by adding some dynamics to the fluctuations of the counting
function and define, for each $\alpha \in (0,2], c>0, a \in
\mathbb{R}_+$, the process
\begin{equation*}
Z_t^{\alpha,c,a} (s) := \frac{N_t^{\alpha,c,a}[0,s t^{1/\alpha}] - s
  t^{1/\alpha}/a}{t^{1/2 \alpha}}, \;\;s \in [0,\infty).
\end{equation*}

\subsection{The covariance structure}
We begin to characterize these processes by identifying the covariance structure.
\begin{lem} \label{covariancelemma}
$\{Z_t^{\alpha,c,a}(s);s\in[0,\infty)\}$ has covariance structure
\begin{align*}
\textrm{Cov}\big(Z_t^{\alpha,c,a}(r),Z_t^{\alpha,c,a}(s)\big)=
\frac{1}{2}\big( f^{\alpha,c,a}(s) + f^{\alpha,c,a} (r) -
f^{\alpha,c,a}(|r-s|) \big).
\end{align*}
\end{lem}
\begin{proof}
By construction
\begin{align*}
N_t^{\alpha,c,a}[0,(r \vee s) t^{1/\alpha}]-N_t^{\alpha,c,a}[0,(r
\wedge s) t^{1/\alpha}]\buildrel {\textrm{dist}} \over {=}
N_t^{\alpha,c,a}[0,|r-s| t^{1/\alpha}].
\end{align*}
Hence, from the definition of $Z_t^{\alpha,c,a}$
\begin{align*}
Z_t^{\alpha,c,a}(|r-s|) \buildrel{\textrm{dist}} \over {=}
Z_t^{\alpha,c,a}(r \vee s)- Z_t^{\alpha,c,a}(r \wedge s),
\end{align*}
which implies that
\begin{align*}
\textrm{Var}\big(Z_t^{\alpha,c,a}(|r-s|)\big) = &
\textrm{Var}\big(Z_t^{\alpha,c,a}(r \wedge s)\big) +
\textrm{Var}\big(Z_t^{\alpha,c,a}(r \vee s)\big) \\ & \quad -
2\textrm{Cov}\big(Z_t^{\alpha,c,a}(r \wedge s), Z_t^{\alpha,c,a}(r
\vee s)\big).
\end{align*}
Rearranging gives
\begin{align*}
&\textrm{Cov}\big(Z_t^{\alpha,c,a}(s), Z_t^{\alpha,c,a}(r) \big) \\ &
\quad = \frac{1}{2}\Big(\textrm{Var}\big(Z_t^{\alpha,c,a}(r \wedge
s)\big) + \textrm{Var}\big(Z_t^{\alpha,c,a}(r \vee s)\big)-
\textrm{Var}\big(Z_t^{\alpha,c,a}(|r-s|)\big)\Big)
\\ & \quad = \frac{1}{2\,t^{1/\alpha}}\Big(V_t^{\alpha,c,a}[s t^{1/\alpha}]+V_t^{\alpha,c,a}[
rt^{1/\alpha}]-V_t^{\alpha,c,a}[ |r - s|t^{1/\alpha}]\Big).
\end{align*}
On referring back to the definition of $f^{\alpha,c,a}(\cdot)$ we see
that this last statement is equivalent to the result of the Lemma.
\end{proof}

Note that the covariance does not depend on t.

\subsection{Convergence of finite dimensional distributions}
Given the covariance structure of $Z_t^{\alpha,c,a}$ and the
identification of its Gaussian one dimensional marginal distributions
the natural next step is to consider the finite dimensional
distributions.

\begin{prop}\label{finitedimprop}
Let $\{G^{\alpha,c,a}(s) : s \in [0,\infty) \}$ be a centered Gaussian
process with covariance structure
\begin{eqnarray}
\textrm{Cov}\big(G^{\alpha,c,a}(s_i),G^{\alpha,c,a}(s_j)\big)=
\frac{1}{2}\big( f^{\alpha,c,a}(s_i) + f^{\alpha,c,a} (s_j) -
f^{\alpha,c,a}(|s_i-s_j|) \big). \label{covfunctiong}
\end{eqnarray}
For any $0 \leq s_1 \leq s_2 \leq \cdots \leq s_n < \infty $ we have
\begin{equation}
(Z_t^{\alpha,c,a}(s_1), Z_t^{\alpha,c,a}(s_2), \ldots,
Z_t^{\alpha,c,a}(s_n)) \Rightarrow (G^{\alpha,c,a}(s_1),
G^{\alpha,c,a}(s_2), \ldots, G^{\alpha,c,a}(s_n))\nonumber
\end{equation}
as $t \rightarrow \infty$.
\end{prop}
\begin{proof}
As previously noted, the mean and covariance structure of
$Z_t^{\alpha,c,a}(s)$ are not dependent on t. Therefore, all that
remains is to show that, in the limit as $t \rightarrow \infty$, the
joint distributions are Gaussian. We again make use of the cumulants.

Recall that given a random vector ${\bf Y} := (Y_1,Y_2, \ldots, Y_n)
\in \mathbb{R}^n$, the joint cumulants of ${\bf Y}$ denoted
$C_{m_1,m_2, \ldots, m_n}({\bf Y})$ are defined via the $m_j$'th
partial derivatives of the logarithm of the characteristic function of
${\bf Y}$. That is,
\begin{eqnarray}
C_{m_1,m_2, \ldots, m_n}\!({\bf Y}) := \!\!
\Bigg(\!\frac{\partial}{\partial \big( i \theta_1\big)}
\!\Bigg)^{m_1}\!\!\!  \Bigg(\!\frac{\partial}{\partial \big(i
  \theta_2\big)}\! \Bigg)^{m_2}\!\!\!\!\! \!\! \cdots \!\!
\Bigg(\!\frac{\partial}{\partial \big( i \theta_n
\big)}\!\!\!\Bigg)^{m_n}\!\!\!\! \log \mathbb{E}\! \Big[\!\exp \big(\!
\sum_{j=1}^n \!\! i \theta_j Y_j \big)\!\Big]\! \Big|_{\bf{\theta}=0}.
\nonumber
\end{eqnarray}
If
\begin{eqnarray}
C_{0,0,\ldots, \underbrace{1}_{\textrm{i'th}},\ldots,0} ({\bf Y}) &=&
\mathbb{E}[Y_i]=0 \nonumber \\
C_{0,0,\ldots, \underbrace{2}_{\textrm{i'th}},\ldots,0} ({\bf Y})&=&
\textrm{Var}[Y_i] = \Sigma_{ii} \nonumber \\
C_{0,\ldots, \underbrace{1}_{\textrm{i'th}}, \ldots, 0,\ldots,
  \underbrace{1}_{\textrm{j'th}},\ldots,0} ({\bf Y}) &=&
\textrm{Cov}[Y_i,Y_j] = \Sigma_{ij} \nonumber
\end{eqnarray}
and in particular
\begin{eqnarray}
C_{m_1,m_2,\ldots,m_n} ({\bf Y}) = 0, \qquad \textrm{whenever} \quad
\sum_{i=1}^n m_i \geq 3 ,\nonumber
\end{eqnarray}
then ${\bf Y}$ has a multivariate normal$(0 , \Sigma)$ distribution. To
prove the Proposition it is enough to show that
\begin{eqnarray}
(Z_t^{\alpha,c,a}(s_1), Z_t^{\alpha,c,a}(s_2)) \rightarrow
\textrm{MultivariateNormal}(0,\Sigma^{\alpha,c,a}) \label{claim1}
\end{eqnarray}
in distribution as $t \rightarrow \infty$, where $\Sigma^{\alpha,c,a}$
is the $2 \times 2$ covariance matrix
\begin{align*}
\begin{pmatrix}\!\!\!\! f^{\alpha,c,a}(s_1)\!\!\!\!\!
 &\!\!\!\!\! \frac{1}{2}\!\big(\! f^{\alpha,c,a}\!(\!s_1\!) \!+\!
 f^{\alpha,c,a} \!(\!s_2\!) \!-\! f^{\alpha,c,a}\!(\!|s_1-s_2|\!)\! \big) \\
\frac{1}{2}\!\big(\! f^{\alpha,c,a}\!(\!s_1\!)\! + \!f^{\alpha,c,a}\!
(\!s_2\!)\! -\! f^{\alpha,c,a}\!(\!|s_1-s_2|\!)\! \big)\! & \!f^{\alpha,c,a}(s_2)\!\!\!\!
\end{pmatrix}.
\end{align*}

We begin by computing the characteristic function of
\begin{align*}
\big( N_t^{\alpha,c,a}[0,s_1 t^{1/\alpha}], N_t^{\alpha,c,a}[0,s_2 t^{1/\alpha}]\big).
\end{align*}
Using the independence of the individual particles we have
\begin{align*}
  &\mathbb{E}_{\mathbb{P}}\Big[\exp(i \theta_1 N_t^{\alpha,c,a}[0 ,
  s_1 t^{1/\alpha}] + i \theta_2 N_t^{\alpha,c,a}[0 , s_2 t^{1/\alpha}])\Big]\\
  &= \mathbb{E}_{\mathbb{P}}\Big[\exp\,i \Big((\theta_1 + \theta_2)
  N_t^{\alpha,c,a}[0 , s_1 t^{1/\alpha}] + \theta_2
  N_t^{\alpha,c,a}[s_1 t^{1/\alpha} , s_2 t^{1/\alpha}]\Big)\Big] \\
  &= \!\!\! \prod_{j=-\infty}^{\infty}\!\!\Big(
  e^{i(\theta_1+\theta_2)}
  \mathbb{P}\big[X_j^{\alpha,c}\!(t)\!+\!u_j\! \in [0 , s_1
  t^{1/\alpha}]\!\big] \\ & \quad \quad + e^{i \theta_2}
  \mathbb{P}\big[X_j^{\alpha,c}(t)+u_j \in [s_1 t^{1/\alpha} , s_2
  t^{1/\alpha}]\big] + \mathbb{P}\big[X_j^{\alpha,c}\!(t)\!+\!u_j
  \!\notin [0 , s_2 t^{1/\alpha}]\big]\Big).
\end{align*}
For ease of notation we will henceforth let
\begin{align*}
q_j(s_l,s_r) := \mathbb{P}[X_j^{\alpha,c}(t)+u_j \in [s_l t^{1/\alpha}
, s_r t^{1/\alpha}]], \qquad 0\leq s_l \leq s_r < \infty.
\end{align*}
The joint cumulants are given by
\begin{align*}
  &C_{m_1,m_2}\big(N_t^{\alpha,c,a}[0 , s_1 t^{1/ \alpha}],
  N_t^{\alpha,c,a}[0 , s_2 t^{1/ \alpha}]\big) \\ & \!=
  \sum_{j=\infty}^{\infty}\!\!\! \Big(\frac{\partial}{\partial (i
    \theta_1)} \Big)^{m_1} \!\! \Big(\frac{\partial}{\partial
    (i\theta_2)} \Big)^{m_2}\!\!\! \log\! \Big(e^{i (\theta_1
    +\theta_2)}\!q_j\!(0,s_1\!)\! +\! e^{i \theta_2}\!
  q_j\!(s_1,s_2\!)\! + \!1-q_j\!(0,s_2\!)\! \Big)\!\Big|_{{\bf
      \theta}=0}.
\end{align*}
Now using the fact that
\begin{align*}
  &\Big(\frac{\partial}{\partial (i \theta_1)} \Big)^{m_1} \!\!
  \Big(\frac{\partial}{\partial (i\theta_2)} \Big)^{m_2} e^{i
    (\theta_1 +\theta_2)}q_j(0,s_1) + e^{i \theta_2} q_j(s_1,s_2) +
  1-q_j(0,s_2) \\ =\,\, & \quad e^{i(\theta_1+\theta_2)}q_j(0,s_1)
  \qquad \forall\,\, m_1,m_2 \, \textrm{ s.t. }m_1 \geq 1
\end{align*}
and
\begin{align*}
  &\Big(\frac{\partial}{\partial (i\theta_2)} \Big)^{m_2} e^{i
    (\theta_1 +\theta_2)}q_j(0,s_1) + e^{i \theta_2} q_j(s_1,s_2) +
  1-q_j(0,s_2) \\ =\,\, & \quad e^{i(\theta_1+\theta_2)}q_j(0,s_1) +
  e^{i \theta_2}q_j(s_1,s_2) \qquad \forall \,\,m_2,
\end{align*}
along with
\begin{align*}
\Big(e^{i (\theta_1 +\theta_2)}q_j(0,s_1) + e^{i \theta_2}
q_j(s_1,s_2) + 1-q_j(0,s_2)\Big) \Big|_{{\bf \theta}=0} = 1,
\end{align*}
we deduce the following generalization of the recursive relation
\eqref{recursive}, with obvious short-hand notation for the joint
cumulants
\begin{align*}
  C_{m_1,m_2} &= \sum_{\substack{k,l: \\ 2\leq k+l \leq m_1+m_2 -1}}
  \!\!\!\!\! \beta_{k,l,m_1,m_2}\, C_{k,l}\\ & \,\, + \,\,
  (-1)^{m_1+m_2}(m_1+m_2 -1)! \sum_{j=-\infty}^{\infty}\!
  q_j\!(0,s_1)\! - \!q_j\!(0,s_1)^{m_1}\!q_j(0,s_2)^{m_2}.
\end{align*}
Now suppose that $C_{k,l} = o(t^{\frac{k+l}{2 \alpha}})$ for all $k,l$
such that $k+l \in \{3,4,\ldots,m_1+m_2-1\}$ and without loss of
generality assume that $m_1+m_2$ is even. Since
\begin{equation*}
0 \leq q_j(0,s_1) - q_j(0,s_1)^{m_1}q_j(0,s_2)^{m_2} \leq q_j(0,s_1),
\end{equation*}
if we assume that the above induction hypothesis holds, then we have
\begin{equation}
o(t^{\frac{m_1+m_2-1}{2 \alpha}})\leq C_{m_1,m_2} \leq
o(t^{\frac{m_1+m_2-1}{2 \alpha}}) + (m_1+m_2-1)!
\frac{s_1t^{1/\alpha}}{a}, \label{recursive2}
\end{equation}
which implies that $C_{m_1,m_2} = o(t^{\frac{m_1+m_2}{2 \alpha}})$
also. We check the third order joint cumulants directly and deduce
that
\begin{align*}
o(t^{\frac{3}{2 \alpha}}) & \leq \,\,C_{k,l} \leq\,\, o(t^{\frac{3}{2
    \alpha}}) + 2\Big(\frac{s_2 t^{1/\alpha}}{a}\Big) \quad
\textrm{when } \, k+l =3,
\end{align*}
since the variances and covariance of $(N_t^{\alpha,c,a}[0,s_1
t^{1/\alpha}], N_t^{\alpha,c,a}[0,s_2 t^{1/\alpha}]$) (i.e.
$C_{2,0},C_{0,2},C_{1,1}$) grow at most like $t^{1/\alpha}$ as $t
\rightarrow \infty$.  Therefore, by induction, whenever $m_1 + m_2
\geq 3$ we have
\begin{eqnarray}
\frac{C_{m_1,m_2}(N_t^{\alpha,c,a}[0 , s_1 t^{1/ \alpha}],
  N_t^{\alpha,c,a}[0 , s_2 t^{1/ \alpha}])}{t^{(m_1+m_2)/ 2 \alpha}}
\rightarrow 0 \quad \textrm{as } t \rightarrow \infty. \nonumber
\end{eqnarray}
In terms of the joint cumulants of $Z_t^{\alpha,c,a}$ this implies
\begin{eqnarray}
  C_{m_1,m_2}(Z_t^{\alpha,c,a}(s_1),Z_t^{\alpha,c,a}(s_2)) \rightarrow
  0 \quad \textrm{as} \,\, t \rightarrow \infty \quad
  \textrm{whenever} \,\, m_1 + m_2 \geq 3,
  \nonumber
\end{eqnarray}
from which the claim \eqref{claim1} and the Proposition follow.
\end{proof}

\subsection{A functional limit for the fluctuation process}

In order to give a functional limit result we consider a continuous
approximation $\{\widehat{Z}^{\alpha,c,a}_t(s): s \in [0,\infty)\}$ to
the process $\{Z_t^{\alpha,c,a}(s): s \in [0,\infty)\}$. Let
\begin{eqnarray}
\widehat{Z}_t^{\alpha,c,a} (s) := \frac{\widehat{N}_t^{\alpha,c,a}[0,s
  t^{1/\alpha}] - s t^{1/\alpha}/a}{t^{1/2 \alpha}},
\end{eqnarray}
where $\widehat{N}_t^{\alpha,c,a}[0,s t^{1/\alpha}]$ is defined to be
equal to $N_t^{\alpha,c,a}[0,s t^{1/\alpha}]$ except at the points of
discontinuity where we linearly interpolate. Let $C[0,1]$ be the space
of continuous real valued functions on $[0,1]$ equipped with the
uniform topology.  We shall denote the measure induced by
$\{\widehat{Z}_t^{\alpha,c,a}(s) : s \in [0,1] \}$ on the space
$(C[0,1], \mathcal{B}(C[0,1]))$ by $Q_t^{\alpha,c,a}$. To simplify
notation we restrict attention to the interval $[0,1]$ but note that
the ensuing functional limit theorem extends trivially to any finite
real indexing set.  The remainder of this section is devoted to
establishing the following weak convergence result.
\begin{thm}\label{funclim}
  Let $Q^{\alpha,c,a}$ be the law of the centered Gaussian process
  $\{G^{\alpha,c,a}(s) : s \in[0,1]\}$ introduced in the statement of
  Proposition \ref{finitedimprop}.  Then
\begin{eqnarray}
Q_t^{\alpha,c,a} \Rightarrow Q^{\alpha,c,a} \nonumber
\end{eqnarray}
as $t \rightarrow \infty$.
\end{thm}
\begin{proof}
Note that by definition
\begin{eqnarray}
|\widehat{Z}_t^{\alpha,c,a}(s)- Z_t^{\alpha,c,a}(s)| \leq \frac{1}{t^{1/2\alpha}}. \label{lessthant}
\end{eqnarray}
Thus, as $t \rightarrow \infty$, the finite dimensional distributions
of $\widehat{Z}_t^{\alpha,c,a}(s)$ must converge to the finite
dimensional distributions of the limiting process $G^{\alpha,c,a}(s)$
to which those of $Z^{\alpha,c,a}_t$ converge. Hence, immediately from
Proposition \ref{finitedimprop} we have that for any $0 \leq s_1 \leq
s_2 \leq \cdots \leq s_n < \infty$,
\begin{equation*}
(\widehat{Z}_t^{\alpha,c,a}(s_1), \widehat{Z}_t^{\alpha,c,a}(s_2),
\ldots, \widehat{Z}_t^{\alpha,c,a}(s_n)) \Rightarrow
(G^{\alpha,c,a}(s_1), G^{\alpha,c,a}(s_2), \ldots,
G^{\alpha,c,a}(s_n))
\end{equation*}
as $t \rightarrow \infty$.  Therefore, by a well known result of
Prohorov's (see e.g. \cite{billings}), the proposed functional limit
theorem holds if the sequence of measures $\{Q_t^{\alpha,c,a}\}$ is
tight. Indeed this tightness requirement follows from Proposition
\ref{tightness} given below.
\end{proof}

The sufficient conditions for tightness verified below are stated in
terms of the distributions
\begin{eqnarray}
\mathbb{P}(\widehat{Z}_t^{\alpha,c,a} \in A) = Q_t^{\alpha,c,a}(A)
\qquad A \in \mathcal{B}(C[0,1]),\nonumber
\end{eqnarray}
and the modulus of continuity, which in this case is defined by
\begin{equation*}
w(\widehat{Z}_t^{\alpha,c,a}, \delta) := \sup_{|s-r| \leq \delta} |
\widehat{Z}_t^{\alpha,c,a}(s)-\widehat{Z}_t^{\alpha,c,a}(r) |, \qquad
\delta \in (0,1].
\end{equation*}

\begin{prop}\label{tightness}
Given $\, \epsilon,\, \lambda >0 \,\, \exists \delta > 0, \,\, t_0 \in \mathbb{N}$ such that
\begin{eqnarray}
\mathbb{P}[w(\widehat{Z}_t^{\alpha,c,a}, \delta) \geq \lambda ] \quad \leq \quad
\epsilon, \qquad \textrm{for }\, t \geq t_0. \nonumber
\end{eqnarray}
\end{prop}
\noindent{Proposition \ref{tightness} is proven via the following
  series of Lemmas.}
\begin{lem} \label{lemtwo}
Suppose $0\leq u \leq r \leq s \leq v \leq 1$, then
\begin{eqnarray}
|\widehat{Z}_t^{\alpha,c,a}(s) - \widehat{Z}_t^{\alpha,c,a}(r)| \leq
|\widehat{Z}_t^{\alpha,c,a}(v) - \widehat{Z}_t^{\alpha,c,a}(u)| + (v-u)t^{1/2\alpha}. \nonumber
\end{eqnarray}
\end{lem}
\begin{proof}
Clearly, by construction
\begin{eqnarray}
0 \leq \widehat{N}_t^{\alpha,c,a}[0, s t^{1/\alpha}]-
\widehat{N}_t^{\alpha,c,a}[0, r t^{1/\alpha}] \leq
\widehat{N}_t^{\alpha,c,a}[0, v t^{1/\alpha}]-
\widehat{N}_t^{\alpha,c,a}[0, u t^{1/\alpha}]. \nonumber
\end{eqnarray}
Therefore, using the definition of $\widehat{Z}^{\alpha,c,a}_t$, we have
\begin{eqnarray}
0 \leq \widehat{Z}_t^{\alpha,c,a}\!(s) -
\widehat{Z}_t^{\alpha,c,a}\!(r) + \frac{(s-r)}{a}t^{1/2\alpha}
 \leq \widehat{Z}_t^{\alpha,c,a}\!(v) -
 \widehat{Z}_t^{\alpha,c,a}\!(u) + \frac{(v-u)}{a}t^{1/2\alpha}.
 \nonumber
\end{eqnarray}
The result follows by rearranging, using the facts $a \in \mathbb{R}_+$, $(v-u) \geq (s-r) \geq 0$  and then
considering separately each case
\begin{align*}
\widehat{Z}_t^{\alpha,c,a}(s) - \widehat{Z}_t^{\alpha,c,a}(r) \geq 0 \qquad \textrm{or} \qquad
\widehat{Z}_t^{\alpha,c,a}(s) - \widehat{Z}_t^{\alpha,c,a}(r) <0.
\end{align*}
\end{proof}
\begin{lem} \label{lemone}
\begin{eqnarray}
|\widehat{Z}_t^{\alpha,c,a}(s) - \widehat{Z}_t^{\alpha,c,a}(r)| \leq
\frac{2}{t^{1/ 2 \alpha}} + |Z_t^{\alpha,c,a}(s) -
Z_t^{\alpha,c,a}(r)|.
\end{eqnarray}
\end{lem}
\begin{proof}
Follows from \eqref{lessthant} and an application of the triangle inequality.
\end{proof}

To obtain results on the distribution of the modulus of continuity for
our sequence of processes $\{\widehat{Z}_t^{\alpha,c,a}\}$ we divide
the interval $[0,1]$ into $m$ disjoint subintervals of length
approximately $\delta$ as follows. Let
\begin{eqnarray}
0= r_0 < r_{k_1} < \cdots < r_{k_{m-1}} < r_{k_m}=1, \label{partition}
\end{eqnarray}
where we define
\begin{eqnarray}
r_i &:=& \frac{i}{t^{1/\alpha}}, \qquad i \in \mathbb{N} \nonumber \\
k_j &:=& j \lceil \delta t^{1/\alpha} \rceil, \qquad j\in \{0,1,2, \ldots,m-1\} \nonumber
\end{eqnarray}
and $\lceil \cdot \rceil$ denotes the ceiling function.
We have
\begin{eqnarray}
\delta \leq r_{k_j} - r_{k_{j-1}} \leq \delta + \frac{1}{t^{1/\alpha}}
\qquad j \in \{1,2, \ldots , m-1 \}, \nonumber
\end{eqnarray}
with the subintervals $[r_{i-1}, r_{i}]$, $i \in \mathbb{N}$ typically
being shorter in length.
\begin{lem}\label{lem4}
\begin{equation*}
\mathbb{P}\Big[ w(\widehat{Z}^{\alpha,c,a}_t , \delta) \geq
\lambda\Big] \leq \sum_{j=0}^{m-1}\! \mathbb{P}\Big[\!
\max\limits_{k_j \leq i \leq k_{j+1}}\!\!
|Z_t^{\alpha,c,a}(r_i)-Z_t^{\alpha,c,a}(r_{k_j})|
\geq \frac{\lambda}{9} - \frac{7}{3\,t^{1/ 2 \alpha}} \Big].
\end{equation*}
\end{lem}
\begin{proof}
Given the partition \eqref{partition}, standard methods (see Theorem
7.4 of \cite{billings}) yield
\begin{equation}
\mathbb{P}\Big[w(\widehat{Z}_t^{\alpha,c,a}, \delta) \geq \lambda
\Big] \leq \sum_{j=0}^{m-1} \mathbb{P} \Big[\sup\limits_{r_{k_j} \leq
s \leq r_{k_{j+1}}} |
\widehat{Z}_t^{\alpha,c,a}(s)-\widehat{Z}_t^{\alpha,c,a}(r_{k_j})|
\geq \frac{\lambda}{3} \Big]. \label{billing}
\end{equation}
By the triangle inequality we have
\begin{equation}
|\widehat{Z}_t^{\alpha,c,a}(s)-\widehat{Z}_t^{\alpha,c,a}(r_{k_j})|
\leq |\widehat{Z}_t^{\alpha,c,a}(s)-\widehat{Z}_t^{\alpha,c,a}(r_{i})|
+
|\widehat{Z}_t^{\alpha,c,a}(r_i)-\widehat{Z}_t^{\alpha,c,a}(r_{k_j})|.
\label{triangle}
\end{equation}
Now if $s \in [r_{k_j}, r_{k_{j+1}}]$, then either
\begin{equation*}
\widehat{Z}_t^{\alpha,c,a}(s)=\widehat{Z}_t^{\alpha,c,a}(r_{i}) \quad
\textrm{for some} \quad i \in \mathbb{N}
\end{equation*}
immediately simplifying \eqref{triangle}, or $\exists i
\in \mathbb{N}$ such that
\begin{equation*}
r_{k_j} \leq r_{i-1} < s < r_i \leq r_{k_{j+1}}
\end{equation*}
in which case from Lemma \ref{lemtwo} we have
\begin{align}
&|\widehat{Z}_t^{\alpha,c,a}(s)-\widehat{Z}_t^{\alpha,c,a}(r_{i})|
\nonumber \\ & \leq
|\widehat{Z}_t^{\alpha,c,a}(r_i)-\widehat{Z}_t^{\alpha,c,a}(r_{i-1})|
+ (r_i - r_{i-1}) t^{1/ 2 \alpha} \nonumber \\
& =
|\widehat{Z}_t^{\alpha,c,a}(r_i)-\widehat{Z}_t^{\alpha,c,a}(r_{i-1})|
+\frac{1}{t^{1/2\alpha}} \nonumber \\
& \leq |\widehat{Z}_t^{\alpha,c,a}(r_i)-\widehat{Z}_t^{\alpha,c,a}(r_{k_j})|
+ |\widehat{Z}_t^{\alpha,c,a}(r_{i-1}) -
\widehat{Z}_t^{\alpha,c,a}(r_{k_j})| +\frac{1}{t^{1/2\alpha}}.
\label{inbetween}
\end{align}
Therefore, using the inequality \eqref{triangle} in conjunction with
\eqref{inbetween} and Lemma \ref{lemone}, we have that for $s, r_i \in
[r_{k_j}, r_{k_{j+1}}]$
\begin{align*}
& \Big|\widehat{Z}_t^{\alpha,c,a}(s)-\widehat{Z}_t^{\alpha,c,a}(r_{k_j})\Big| \\
\leq \quad &
\Big|\widehat{Z}_t^{\alpha,c,a}(r_i)-\widehat{Z}_t^{\alpha,c,a}(r_{k_j})\Big|
+
\Big|\widehat{Z}_t^{\alpha,c,a}(r_{i-1})-\widehat{Z}_t^{\alpha,c,a}(r_{k_j})\Big|
+\frac{1}{t^{1/2\alpha}}
\\ &+
\Big|\widehat{Z}_t^{\alpha,c,a}(r_i)-\widehat{Z}_t^{\alpha,c,a}(r_{k_j})\Big|
\\ \leq \quad &
3 \max\limits_{k_j \leq i \leq k_{j+1}} \Big|
\widehat{Z}_t^{\alpha,c,a}(r_i)-\widehat{Z}_t^{\alpha,c,a}(r_{k_j})\Big|
+\frac{1}{t^{1/2\alpha}} \\
\leq \quad& 3 \,\Big[\max\limits_{k_j \leq i \leq k_{j+1}}
\Big|Z_t^{\alpha,c,a}(r_i) - Z_t^{\alpha,c,a}(r_{k_j})\Big| +
\frac{7}{3t^{1/2 \alpha}}\Big].
\end{align*}
Thus
\begin{align*}
&\mathbb{P}\Big[\sup\limits_{r_{k_j} \leq s \leq r_{k_{j+1}}}
\Big|\widehat{Z}_t^{\alpha,c,a}(s) -
\widehat{Z}_t^{\alpha,c,a}(r_{k_j})\Big| \geq \frac{\lambda}{3}\Big]
\\ \leq \,\, &\mathbb{P}\Big[ \max\limits_{k_j \leq i \leq k_{j+1}}
\Big|Z_t^{\alpha,c,a}(r_i) - Z_t^{\alpha,c,a}(r_{k_j})\Big| +
\frac{7}{3t^{1/2 \alpha}} \geq \frac{\lambda}{9} \Big].
\end{align*}
Substituting this last inequality into \eqref{billing} gives the
statement of the Lemma.
\end{proof}
Now that we have reduced the study of the distribution of the modulus
of continuity to that of the maximum fluctuation over our constructed
subintervals we can progress by introducing a maximal inequality. In
order to do this we use the following known result
taken from \cite{billings} and paraphrased for use here.
\begin{thm}\label{billingsthm}
Consider a sequence of random variables $\{\xi_q\}_{q \geq 1}$ and the
associated partial sums
\begin{equation*}
S_u := \sum_{q=1}^u \xi_q \qquad S_0 := 0.
\end{equation*}
Let
\begin{equation*}
M_n := \max\limits_{1 \leq u \leq n} |S_u|.
\end{equation*}
If
\begin{equation*}
\mathbb{P}\big[ |S_v - S_u | \geq \gamma \big] \leq \frac{1}{\gamma^2}
\Big( \sum_{u <l \leq v} b_l\Big)^2 \qquad 0 \leq u \leq v \leq n
\end{equation*}
for $\gamma >0$ and some $b_1,b_2,\ldots, b_n \in \mathbb{R}_+$, then
\begin{equation*}
\mathbb{P}[M_n \geq \gamma] \leq \frac{\kappa}{\gamma^2} \Big( \sum_{0
<l \leq n} b_l\Big)^2,
\end{equation*}
where $\kappa$ is a constant.
\end{thm}
\begin{proof}
See \cite{billings}, Theorem 10.2
\end{proof}
\begin{lem}\label{maxlem}
\begin{align*}
\mathbb{P}\big[\max\limits_{1 \leq u \leq n}
|Z_t^{\alpha,c,a}(r_{k_j+u}) - Z_t^{\alpha,c,a}(r_{k_j}) | \geq
\gamma] \leq \frac{\kappa}{a \gamma^2} \big(n t^{-1/2\alpha} \big)^2, \end{align*}
where $\kappa$ is constant and $\gamma>0$.
\end{lem}
\begin{proof}
Let
\begin{equation*}
\xi_q := Z_t^{\alpha,c,a}(r_{k_j+q}) - Z_t^{\alpha,c,a}(r_{k_j +q-1}).
\end{equation*}
Then
\begin{equation*}
S_u := \sum_{q=1}^u \xi_q = Z_t^{\alpha,c,a}(r_{k_j+u}) -
Z_t^{\alpha,c,a}(r_{k_j}) \qquad S_0:=0.
\end{equation*}
In this case, for $v,u \in \mathbb{N}$, applying Chebyshev's
inequality and using the definition of $f^{\alpha,c,a}(\cdot)$ and the
upper bound given at \eqref{markovbound}, we have
\begin{align*}
 \mathbb{P}\big[ |S_v - S_u | \geq \gamma \big]
= \quad & \mathbb{P}\big[ |Z_t^{\alpha,c,a}(r_{k_j +v}) -
Z_t^{\alpha,c,a}(r_{k_j +u})| \geq \gamma\big] \\
\leq \quad & \frac{f^{\alpha,c,a}(|r_{k_j+v} - r_{k_j+u}|)}{\gamma^2}\\
\leq
\quad & \frac{(v-u)t^{-1/\alpha}}{a \gamma^2} \\
\leq \quad & \frac{\big[ (v-u) t^{-1/2\alpha} \big]^2}{a \gamma^2}.
\end{align*}
Thus we can take $b_l = \frac{t^{-1/2 \alpha}}{a}$ for $l = 1,2,
\ldots, n$ and apply Theorem \ref{billingsthm} which gives the maximal
inequality of the Lemma.
\end{proof}
Concluding the proof of Proposition \ref{tightness} is now straightforward.
\\[0.2cm]
{\it Proof of Proposition \ref{tightness}.}
Taking $n = k_{j+1}-k_j$ in the statement of Lemma \ref{maxlem}
gives
\begin{align*}
\mathbb{P}\big[\max\limits_{k_j \leq i \leq k_{j+1}} | Z_t(r_i) -
Z_t(r_{k_j})| \geq  \gamma \big]  \leq  \quad &\frac{\kappa}{a
\gamma^2} \big( r_{k_{j+1}}-r_{k_j}\big)^2 \\
 \leq \quad &\frac{\kappa}{a \gamma^2} \big( \delta + t^{-1/\alpha} \big)^2.
\end{align*}
 Substituting this last inequality with $\gamma= \big(
 \frac{\lambda}{9} - \frac{7}{3t^{1/ 2 \alpha}}\big)$ (which is
 strictly positive for sufficiently large $t$) into the inequality
 given by Lemma \ref{lem4} gives
 \begin{align*}
 \mathbb{P}\big[w(Z_t^{\alpha,c,a}, \delta) \geq \lambda \big] &\leq
 \frac{81 m \kappa}{a} \frac{\big(\delta t^{1/2\alpha} + t^{-1/2
 \alpha}\big)^2}{\big( \lambda t^{1/2\alpha} - 21 \big)^2}.
 \end{align*}
On solving the appropriate quadratic equation we find that we can make
this upper bound less than $\epsilon$ by choosing $\delta$ from
\begin{equation*}
\Big(-t^{-1/2\alpha} - c_{\epsilon}\big(\lambda - 21
t^{-1/2\alpha}\big)\, , \, -t^{-1/2\alpha} + c_{\epsilon}\big(\lambda
- 21 t^{-1/2\alpha}\big)\Big)  \bigcap \Big(0,1\Big),
\end{equation*}
where $c_{\epsilon} = \sqrt{\frac{a \epsilon}{81 m \kappa}}$. Since
this interval is non-empty for sufficiently large t this completes the
proof. \hspace{\stretch{1}}$\square$

\subsection{Properties of the limiting process $G^{\alpha,c,a}(s)$}

We have constructed a family $\{(G^{\alpha,c,a}(s),\, s
\in[0,\infty)), \alpha \in(0,2], c>0, a \in \mathbb{R}_+ \}$ of
centered, continuous, real-valued Gaussian processes with the
inherited covariance structure
\begin{align}
Cov(G^{\alpha,c,a}(s),G^{\alpha,c,a}(r)) = \frac{1}{2}\big(
f^{\alpha,c,a}(s) + f^{\alpha,c,a} (r) - f^{\alpha,c,a}(|s-r|)
\label{covfunctiong2}
\big).
\end{align}
It is clear that the processes are recurrent for $1<\alpha\leq 2$ as
the number variance saturates and the stationary distribution is
normal with mean zero and variance $\kappa_{sat}(\alpha,c,a,1)$.
We are able to deduce further properties of these limit processes by
using our earlier results on the number variance for the systems of
symmetric $\alpha$-stable processes from which they are constructed.

\begin{prop}
$G^{\alpha,c,a}(s)$ is non-Markovian.
\end{prop}
\begin{proof}
Recall (see e.g. \cite{kallenberg}, Chapter 13) that a Gaussian
process with indexing set $T \subset \mathbb{R}$ and covariance
function $\rho: T \mapsto \mathbb{R}$ is Markov if and only if
\begin{align*}
\rho(s,u) = \frac{\rho(s,r)\rho(r,u)}{\rho(r,r)} \quad \forall\,\, s,u,r \in T.
\end{align*}
It is clear that this relationship does not hold in general for the
covariance function \eqref{covfunctiong2}.
\end{proof}

From the results of previous sections, since $f^{\alpha,c,a}(s) =
V^{\alpha,c,a}_1[s]$, we know that
\begin{align*}
f^{\alpha,c,a}(s) \rightarrow (1-e^{-2c}) \frac{s}{a} \quad \textrm{as
} \alpha \rightarrow 0
\end{align*}
and
\begin{align*}
f^{\alpha,c,a}(s) \buildrel s \rightarrow 0 \over \sim k_{\alpha,c,a}\, s.
\end{align*}
Therefore, $G^{\alpha,c,a}$ appears to start out, for small ``time" as
a scaled Brownian motion and as $\alpha \rightarrow 0 $ this initial
Brownian character prevails for longer. We capture this more precisely
in the following easily verified proposition.

\begin{prop}
1. $\{G^{\alpha,c,a}(s):s\in [0,1]\}$ converges weakly to a Brownian
motion $\{B(\frac{1-e^{-2c}}{a}s):s\in[0,1]\}$ as $\alpha\to 0$.

2. Let $G_{\epsilon}^{\alpha,c,a}(s) =
\epsilon^{-1/2}G^{\alpha,c,a}(\epsilon s)$. Then
$\{G_{\epsilon}^{\alpha,c,a}(s):s\in [0,1]\}$ converges weakly to a
Brownian motion
$\{B(\frac{s}{a}):s\in[0,1]\}$ as $\epsilon\to 0$.
\end{prop}

\begin{rem}{\rm
The covariance structure of
$G^{\alpha,c,a}$ is similar to that of
a Brownian bridge. Recall that the standard Brownian bridge $(B^{br}(s), s \in
[0,a])$ of length $a$, is a centered Gaussian process with covariance
structure
\begin{align*}
Cov(B^{br}(s),B^{br}(r))= s \wedge r - \frac{sr}{a}
\end{align*}
and arises as a weak limit of many empirical processes. In particular,
it may be obtained from the appropriately scaled counting functions of
a Poisson process on $\mathbb{R}$ (see e.g. \cite{kallenberg}). We can
re-write the covariance \eqref{covfunctiong} in the alternative form
\begin{align*}
Cov(G^{\alpha,c,a}(s), G^{\alpha,c,a}(r)) = \frac{s \wedge r}{a} -
\int_0^{\frac{r}{a}} \int_0^{\frac{s}{a}} p_{2/a^{\alpha}}(y,z) \, \ud
y \, \ud z,
\end{align*}
but we see that a precise match would require
$p_{2/a^{\alpha}}(y,z) =1$. }
\end{rem}


\begin{prop}
The process $G^{\alpha,c,a}$ has stationary increments which are negatively correlated.
\end{prop}
\begin{proof}
It is straightforward to see that the increments have zero mean and
that for any $s,r \in [0,\infty)$
\begin{align*}
\textrm{Var}(G^{\alpha,c,a}(s)-G^{\alpha,c,a}(r)) = f^{\alpha,c,a}(|s-r|).
\end{align*}
In addition, for $u \geq 0$ we have
\begin{align}
  &\textrm{Cov}(G^{\alpha,c,a}(s)-G^{\alpha,c,a}(0),
  G^{\alpha,c,a}(r+s+u)-G^{\alpha,c,a}(s+u)) \nonumber \\= \quad &
  \frac{1}{2}\Big(f^{\alpha,c,a}(s+r+u) - f^{\alpha,c,a}(r+u) -\big[
  f^{\alpha,c,a}(s+u)- f^{\alpha,c,a}(u)\big]\Big).  \label{covinc}
\end{align}
Since $f^{\alpha,c,a}$ is a concave function
\begin{align*}
f^{\alpha,c,a}(s+r+u) - f^{\alpha,c,a}(r+u) \leq
f^{\alpha,c,a}(s+u)- f^{\alpha,c,a}(u),
\end{align*}
for all $s,r,u \in [0,\infty)$ so it follows that
this covariance is non-positive.
\end{proof}

\begin{prop}
$G^{\alpha,c,a}$ is not in general self-similar. For any constant $b
\in \mathbb{R}$ we have the relationship
\begin{align*}
G^{\alpha,c,a}(bs) \buildrel \textrm{dist} \over = b^{1/2} \,
G^{\alpha,\frac{c}{b^{\alpha}},a} (s).
\end{align*}
\end{prop}
\begin{proof}
  Both sides of the proposed equation have zero mean
  and a Gaussian distribution.  It is clear from the expression given for $f^{\alpha,c,a}$ at \eqref{fdefn} that the variances/covariances also agree.
\end{proof}

\begin{prop}\label{longrange}
$G^{\alpha,c,a}$ is a long memory (or long range dependent) process in
the sense that the covariance between increments decays as a power law
as the separation between them is increased. More precisely, for
$\alpha \in (0,2)$ we have
\begin{align*}
Cov\big(G^{\alpha,c,a}(s)-G^{\alpha,c,a}(0),
G^{\alpha,c,a}(r+s+u)-G^{\alpha,c,a}(s+u)\big) \buildrel u \rightarrow
\infty \over \sim k\,u^{-(\alpha +1)},
\end{align*}
where $k$ is a constant depending on $\alpha,c,a,s$ and $r$.
\end{prop}

\begin{proof}
  The covariance in question is expressed in terms of the function
  $f^{\alpha,c,a}$ at \eqref{covinc}. Note from \eqref{asymvar} that
  we already know the asymptotic behaviour of the individual
  components of this expression. Applying l'Hopital's rule twice in
  succession yields the given power law.
\end{proof}

We have already mentioned a similarity between the covariance
structure of $G^{\alpha,c,a}$ and that of Brownian motion. More
generally we can draw parallels between our limiting process and
fractional Brownian motion. Recall (see for example
\cite{mandvanness}) that a fractional Brownian motion $(W_H(s), s \geq
0)$ with Hurst parameter $H \in (0,1)$ is a centered, self-similar
Gaussian process with covariance function
\begin{align} \label{covfrac}
 Cov(W_H(s),W_H(r)) = \frac{1}{2}\Big(
  s^{2H} + r^{2H} - |s-r|^{2H}\Big),
\end{align}
The case $H=\frac{1}{2}$
corresponds to a standard Brownian motion.

Note the resemblance between the form of the covariance functions
\eqref{covfunctiong2} and \eqref{covfrac}.  Heuristically, we can
deduce that $f^{\alpha,c,a}(s)$ may be approximated by a function of
the form $\kappa_{\alpha,c,a} s^{2H_{\alpha,c,a}(s)}$ where
$H_{\alpha,c,a}: [0, \infty) \mapsto [0, \frac{1}{2}]$ is a
monotonically decreasing function with initial value
$H_{\alpha,c,a}(0)=\frac{1}{2}$. Thus loosely speaking
$G^{\alpha,c,a}$ can be viewed as a type of fractional Brownian motion
with time varying Hurst parameter. In particular, the long range dependence property of
Proposition \ref{longrange} may be compared to the analogous statement for fractional Brownian motion:
\begin{align*}
Cov\big(W_H(s)-W_H(0), W_H(r+s+u)-W_H(s+u)\big)\buildrel u \rightarrow
\infty \over \sim k u^{2H-2},
\end{align*}
where $k$ is a constant depending on $H, s$ and $r$.
We make the link between the process $G^{\alpha,c,a}$ and fractional Brownian motion
precise with the following statement.
\begin{prop}
For $\alpha \in (0,1)$ let
\begin{align*}
\widetilde{G}^{\alpha,c,a}_b(s) :=
\frac{G^{\alpha,c,a}(bs)}{\sqrt{b^{1-\alpha}}} , \quad s,b \in
[0,\infty).
\end{align*}
Then
\begin{align*}
\{\widetilde{G}^{\alpha,c,a}_b(s)\,:\, s \in [0,\infty)\} \Rightarrow
\{k_{\alpha,c,a}^{1/2}\,W_{\frac{1-\alpha}{2}}(s) \,:\, s \in
[0,\infty)\} \quad \textrm{as } b \rightarrow \infty
\end{align*}
where $k_{\alpha,c,a} = \frac{4c}{a \pi} \Gamma(\alpha -1) \sin
\big(-\alpha \pi /2\big)$.
\end{prop}
\begin{proof}
From \eqref{fdefn}, by applying a Taylor expansion we deduce that
\begin{align*}
\frac{f^{\alpha,c,a}(bs)}{b^{1-\alpha}} &= s^{1-\alpha} \frac{8c}{a
  \pi} \int_0^{\infty} \frac{\sin^2(u/2)}{u^{2-\alpha}} \ud u \\ &
\quad + \frac{4}{a \pi} \int_0^{\infty}\sin^2(u/2)\sum_{j=2}^{\infty}
\frac{(-1)^{j+1}}{j\,!}(2c)^js^{1-2j\alpha}\frac{u^{j\alpha
    -2}}{b^{j\alpha -\alpha}} \ud u.
\end{align*}
Now by the Dominated Convergence Theorem,
\begin{align*}
&\lim_{b \rightarrow \infty} \frac{4}{a \pi}
\int_0^{\infty}\underbrace{\sin^2(u/2)\sum_{j=2}^{\infty}
  \frac{(-1)^{j+1}}{j\,!}(2c)^js^{1-2j\alpha}\frac{u^{j\alpha
      -2}}{b^{j\alpha -\alpha}}}_{h_b(u)} \ud u \\ & \quad =
\frac{4}{a \pi} \int_0^{\infty}\sin^2(u/2) \sum_{j=2}^{\infty} \lim_{b
  \rightarrow \infty}
\frac{(-1)^{j+1}}{j\,!}(2c)^js^{1-2j\alpha}\frac{u^{j\alpha
    -2}}{b^{j\alpha -\alpha}} \ud u, \\ & \quad = 0
\end{align*}
since, setting
$M(u)=sin^2(u/2)(1-\exp(-2c(u/s)^{\alpha}))/u^2+sin^2(u/2)/u^{2-\alpha}$,
we have a positive integrable function such that $|h_b(u)| \leq M(u)$
for all $b \in \mathbb{R}$. This implies that
\begin{align*}
\lim_{b \rightarrow \infty} \frac{f^{\alpha,c,a}(bs)}{b^{1-\alpha}} =
k_{\alpha,c,a} \, s^{1-\alpha}
\end{align*}
which allows us to conclude
\begin{align*}
Cov\big(\widetilde{G}^{\alpha,c,a}_b(r), \widetilde{G}^{\alpha,c,a}_b(s)\big)&= b^{\alpha -1} Cov\big(G^{\alpha,c,a}(br), G^{\alpha,c,a}(bs)\big)\\
& \buildrel {b \rightarrow \infty} \over \rightarrow  \frac{k_{\alpha,c,a}}{2} \big( s^{1-\alpha} + r^{1-\alpha} + |s-r|^{1-\alpha}\big) \\
&= Cov \big(k_{\alpha,c,a}^{1/2}\,W_{\frac{1-\alpha}{2}}(s), k_{\alpha,c,a}^{1/2}\,W_{\frac{1-\alpha}{2}}(r)\big).
\end{align*}
The processes are Gaussian therefore the convergence of finite
dimensional distributions is implied by the convergence of the
covariance functions and tightness follows easily from, for example,
\cite{kallenberg}~Corollary~16.9 and well known
expressions for the even moments.
\end{proof}


\begin{rem}{\rm
  We mention that various other long-range dependent Gaussian
  processes have recently been found to arise from the fluctuations of
  spatially distributed particle systems, see \cite{goro3} and
  references within. Notably, in this context, the spatial particle
  configurations of infinite systems of symmetric $\alpha$-stable
  processes started from a Poisson process on $\mathbb{R}$ have been
  considered. In these cases, fractional Brownian motion with Hurst
  parameter $H=2-\frac{1}{\alpha}, \, \alpha \in (1,2]$ was obtained
  as a scaling limit of the occupation time process (essentially
  scaling the counting function in time rather than in time and space
  as in this paper).}
\end{rem}

\begin{rem}{\rm
  It seems natural to ask whether, in the same fashion as we created
  $G^{\alpha,c,a}$, similar limiting processes could be constructed
  from Johansson's systems of non-colliding Brownian motions.
  Unfortunately, the formula for the averaged number variance given
  for these processes in \cite{kurt:varsat} does not scale in time and
  space in the same convenient way as $V_t^{\alpha,c,a}[L]$ in this
  case. However, as noted by Johansson, letting $t \rightarrow \infty$
  in his model, one obtains the sine kernel determinantal process,
  from which a limiting Gaussian process (parameterized and scaled in
  a completely different way) was constructed in
  \cite{soshnikov:process}.}
\end{rem}

\end{document}